\documentclass[a4paper, reqno, 12pt]{amsart}
\usepackage{amsmath,amssymb,graphicx,url}
\usepackage{titlesec}
\usepackage{geometry}
\usepackage[most]{tcolorbox}
\usepackage{xcolor}
\usepackage{framed}
\usepackage{tikz}
\usetikzlibrary{shapes.geometric, arrows.meta, positioning}
\usepackage{tabularx}
\usepackage{pdflscape}

\usepackage{array}

\usepackage{tcolorbox}
\usepackage{enumitem}  
\setlist[itemize]{nosep}

\definecolor{low}{rgb}{1,1,1}
\definecolor{med}{rgb}{0.8,0.9,1}
\definecolor{high}{rgb}{0.4,0.7,1}

\usepackage{colortbl}

\tikzstyle{startstop} = [rectangle, rounded corners, minimum width=3.5cm, minimum height=1cm,text centered, draw=black, fill=blue!15]
\tikzstyle{process} = [rectangle, minimum width=3.5cm, minimum height=1cm, text centered, draw=black, fill=orange!20]
\tikzstyle{arrow} = [thick,->,>=stealth]

\titleformat{\section}{\large\bfseries}{\thesection}{1em}{}
\titleformat{\subsection}{\normalsize\bfseries}{\thesubsection}{1em}{}

\theoremstyle{plain}
\newtheorem{theorem}{Theorem}[section]

\theoremstyle{definition}

\newtheorem{remark}[theorem]{Remark}

\title[The VIBE Framework]{The VIBE Framework: A Student-Centered Approach to Teaching Knot Theory in Secondary Mathematics}
\author{Ioannis Diamantis}
\address{Department of Data Analytics and Digitalisation,
Maastricht University, School of Business and Economics,
P.O.Box 616, 6200 MD, Maastricht,
The Netherlands.}
\email{i.diamantis@maastrichtuniversity.nl}



\begin{document}

\begin{abstract}
Knot theory, a visual and intuitive branch of topology, offers a unique opportunity to introduce advanced mathematical thinking in secondary education. Despite its accessibility and cross-disciplinary relevance, it remains largely absent from standard curricula. This paper proposes the {\it VIBE framework}, a student-centered approach, structured around four pedagogical pillars: Visual, Inquiry-based, Braided (collaborative), and Embedded (contextualized) learning. Rooted in constructivist theory, VIBE supports cognitive development, spatial reasoning, and mathematical engagement across diverse learners. We present a sequence of low-threshold, high-ceiling activities designed to develop core topological concepts while fostering creativity and exploration. Through qualitative heatmaps, clustering visualizations, and classroom snapshots, we demonstrate how knot theory can be transformed into a powerful medium for inquiry and interdisciplinary connection. We believe that the VIBE framework provides a structured yet adaptable approach that supports the integration of deep, meaningful mathematical experiences into secondary education.
\end{abstract}

\maketitle

\section{Introduction}

Most high school mathematics curricula avoid higher-level mathematics due to its perceived abstraction, formalism, and limited relevance to students’ everyday experiences \cite{Schmidt2006, DarlingHammond2010}. In response, recent shifts in mathematics education emphasize student-centered approaches, such as inquiry-based learning and project-based learning, that prioritize meaning-making, collaboration, and problem solving over rote memorization \cite{Bruner1966, Hmelo2004, thomas2000}. Among these innovations, visual learning has emerged as a powerful pedagogical tool, particularly in contexts that require spatial reasoning and conceptual insight \cite{Arcavi2003, Presmeg2006, Duval1999}. However, translating these progressive goals into everyday classrooms remains a challenge \cite{devlin}. This is where knot theory offers unique promise: its tactile and visual nature, minimal algebraic prerequisites, and surprising relevance to biology, art, and computer science position it as a highly accessible yet intellectually rich topic for secondary education \cite{adams2004, murasugi}.

Topology is the study of shape, continuity, and spatial relationships \cite{Poincare}. It remains almost entirely absent from secondary education. This paper aims to address the gap between abstract mathematical topics like topology and their absence from secondary education curricula by proposing a student-centered framework for their integration. Within topology lies a particularly accessible and intuitive subfield: \emph{knot theory}. Unlike many abstract branches of mathematics, knot theory is tactile, visual, and deeply embedded in both cultural practices (e.g. weaving \cite{nim}, Celtic art \cite{Gross}) and scientific phenomena (e.g. DNA entanglement \cite{flapan}, and fluid dynamics \cite{ricca}). It is a rich context for developing spatial reasoning and abstraction skills, especially through visual engagement with topological ideas \cite{Arcavi2003, Duval1999}. Moreover, Knot Theory provides a unique opportunity to engage students in authentic mathematical inquiry without requiring formal prerequisites. A concise overview of the historical development of knot theory is provided at the start of Section 2 to contextualize its evolution and educational relevance.

This paper introduces the {\bf VIBE framework} as a practical and theoretically grounded model for teaching knot theory in high school mathematics. The framework integrates constructivist principles with global competencies, supporting student engagement through hands-on, collaborative, and interdisciplinary activities. Inspired during the author’s completion of the University Teaching Qualification (BKO) at Maastricht University, this work combines pedagogical vision with implementable classroom strategies, aiming to enrich both mathematics instruction and professional teaching practice.

\smallbreak 

The {VIBE framework} is structured around four pedagogical pillars: 
\begin{itemize}
    \item \textbf{V}isual Learning, leveraging string, diagrams, and animations to build spatial reasoning.
    \item \textbf{I}nquiry-Based Exploration, guiding students to ask open-ended questions and investigate properties through play.
    \item \textbf{B}raided Activities, embedding collaborative, hands-on problem-solving as a core methodology.
    \item \textbf{E}mbedded Contexts, integrating connections to biology, art, and technology to reinforce interdisciplinary relevance.
\end{itemize}

\smallbreak 

The framework draws on the tradition of constructivism in mathematics education, as advocated by Piaget (1952) \cite{Piaget1952}, Vygotsky (1978) \cite{Vygotsky1978}, and Bruner (1966) \cite{Bruner1966}, which posits that students learn best through active engagement, personal meaning-making, and social interaction. Knot theory naturally supports these principles through its concrete representations, pattern recognition, and open-ended tasks. Through this contribution, we aim to position knot theory not only as a compelling mathematical topic, but also as a catalyst for reimagining what school mathematics can be: rich, connected, and creative.

\smallbreak

The paper is organized as follows: Section 2 introduces basic concepts from knot theory, laying the groundwork for classroom applications. Section 3 discusses the theoretical foundation, focusing on constructivism and the role of visual learning. Section 4 introduces the VIBE framework in detail, explaining each of its four pedagogical pillars and how they interact in the classroom. Section 5 provides a sequence of sample lessons and activities, that can be adapted to a variety of teaching contexts. Section 6 outlines practical strategies for implementation and assessment, emphasizing tools that align with visual and inquiry-based learning. Section 7 discusses the broader benefits and challenges of teaching knot theory at the secondary level, addressing issues such as curriculum alignment, teacher preparedness, and student engagement. Finally, Section 8 concludes with reflections on the potential of knot theory to promote mathematical maturity and outlines directions for future curriculum development and research.

\section{Preliminaries of Knot Theory}

\subsection*{A Brief Historical Perspective}

The history of knot theory is itself an educational narrative. Its origins in the 19th century were not purely mathematical but rooted in the physical sciences. Lord Kelvin proposed that atoms might be knotted tubes of ether, inspiring Peter Guthrie Tait to develop systematic knot tables \cite{burton, adams2004}. Though the theory of the ether was later discarded, this pursuit gave rise to one of the first rigorous explorations of mathematical knots. The field gained greater structure in the 20th century through the work of Kurt Reidemeister \cite{Reid}, whose three moves became the backbone of diagrammatic knot equivalence, and it was revolutionized in the 1980's with the discovery of the Jones polynomial, revealing connections between topology, quantum theory, and algebra \cite{LK}. More recently, knot theory has found concrete applications in biology, especially in modelling the behaviour of knotted DNA molecules \cite{Sum}.

These historical milestones, emerging from curiosity, missteps, and breakthroughs, offer powerful teaching moments. They ``humanize'' mathematics, reminding students that math evolves and is shaped by human imagination and error. Teachers can bring this legacy to life by weaving these narratives into classroom activities, connecting students to a living discipline.

\subsection*{What is a Knot?}

Following the intuitive introductions by Adams in \cite{adams2004, adamswhy}, we avoid formal algebraic language in early activities. In everyday language, a knot is something we tie in rope or string. In mathematics, we define a \textit{knot} more precisely:

\begin{tcolorbox}[title=Definition of Knots and Links]
A \textbf{knot} is a closed loop in three-dimensional space that does not intersect itself. A \textbf{link} is a collection of two or more closed loops in three-dimensional space, which may or may not be entangled.
\end{tcolorbox}

\begin{figure}[ht]
    \centering
    \includegraphics[width=0.5\textwidth]{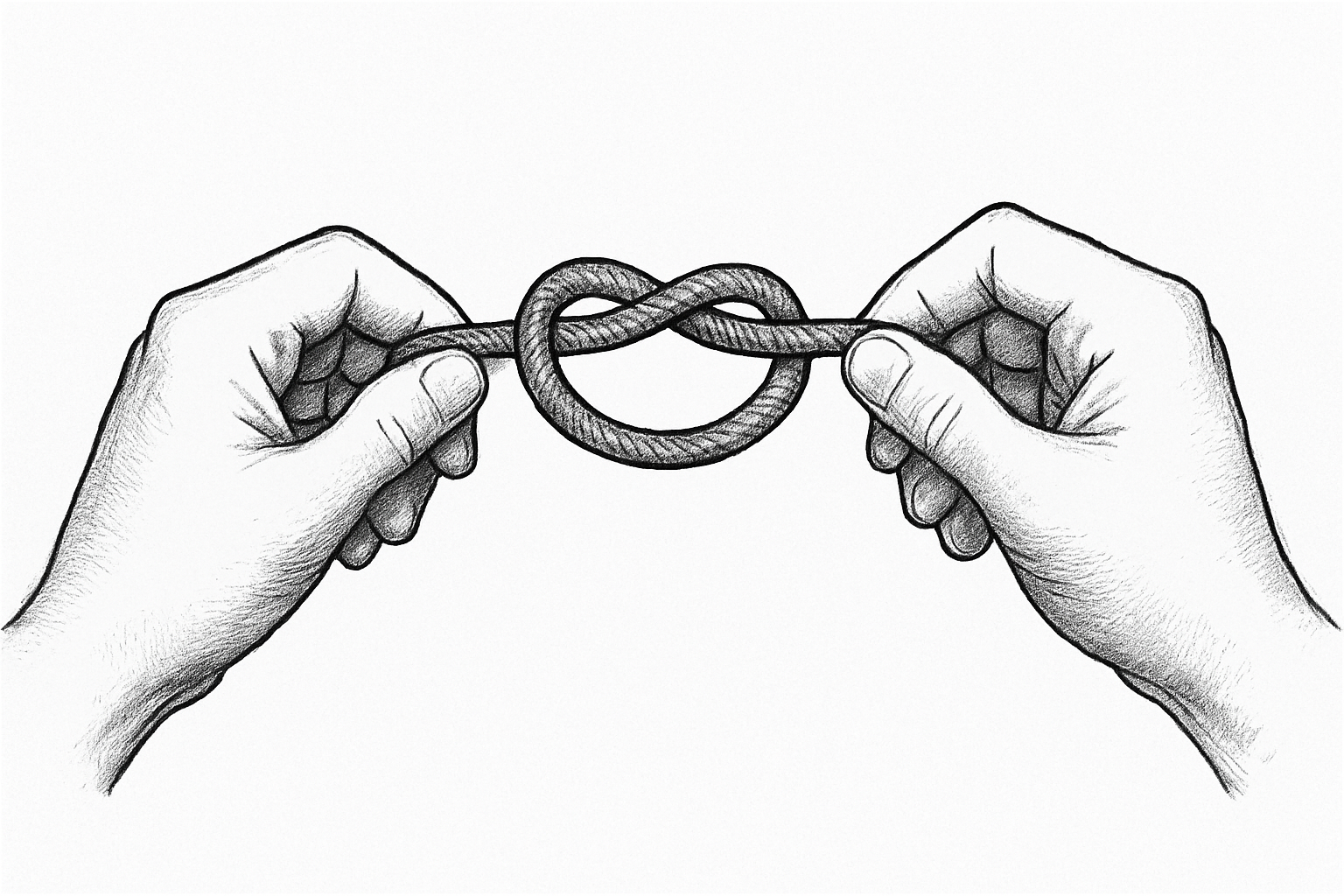}
    \caption{A simple knot formed by tying a loop and connecting the ends.}
    \label{fig:simple-knot}
\end{figure}

Imagine taking a piece of string, tying a knot in it, and then gluing the two ends together so it becomes a continuous loop. The result is a mathematical knot, i.e., no loose ends, just a closed curve that may twist or loop in space (for an illustration see Figure~\ref{fig:simple-knot}).

This definition leaves room for many possibilities: some knots are simple and familiar, while others are more intricate or even involve multiple components. To better understand these distinctions, it's helpful to look at a few basic examples.

\begin{itemize}
    \item The \textbf{unknot} is the simplest possible knot. It’s just a round loop with no twists or crossings.
    \item The \textbf{trefoil knot} is the simplest non-trivial knot. It has three crossings and cannot be deformed into a simple loop without cutting the strand.
    \item A \textbf{link} consists of two or more loops interlaced in space. The simplest such example is the \textbf{Hopf link}, made from two circles linked once.
\end{itemize}

\smallbreak

\begin{figure}[ht]
    \centering
    \includegraphics[width=0.5\textwidth]{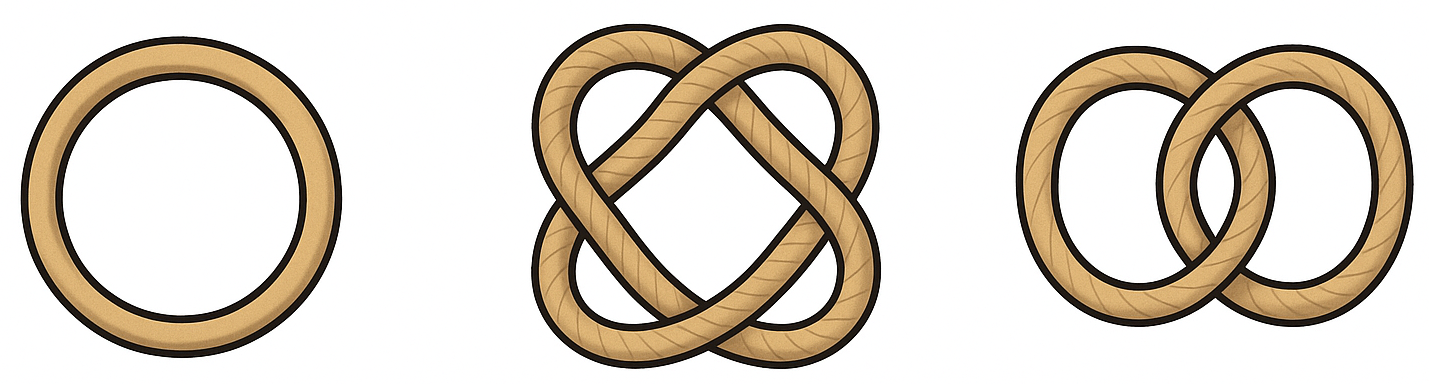}
    \caption{From left to right: the unknot, the Solomon's link, and the Hopf link.}
    \label{fig:basic-knots}
\end{figure}

Figure~\ref{fig:basic-knots} illustrates three fundamental examples in knot theory: the {unknot}, the {Solomon’s link}, a slightly more complex two-component link that demonstrates higher entanglement and the Hopf link.

\subsection*{When Are Two Knots the Same?}

Mathematical knots are considered flexible. We are allowed to stretch, bend, or twist them in space, as long as we don’t cut the string or allow it to pass through itself. These operations may change the appearance of the knot, but not its fundamental structure. This leads to a central question in knot theory: 

\begin{center}
\begin{framed}

{\it When should two knots be considered ``the same''?}
\end{framed}

\end{center}

Before answering this question formally, it helps to build some intuition. Imagine holding a piece of rope: you can twist it, bend it, and move it around in space. But as long as you don’t cut the rope or let it pass through itself, it feels like the same knot. This intuitive idea of deforming one knot into another, without breaking or gluing, captures the topological notion of ``sameness''. Mathematicians have a precise way to describe this kind of transformation.

\begin{tcolorbox}[title=Equivalence and Ambient Isotopy]
Two knots are considered \textbf{equivalent} if one can be ``smoothly'' transformed into the other without cutting or passing strands through each other. This process is called \textbf{ambient isotopy}.
\end{tcolorbox}

\begin{figure}[ht]
    \centering
    \includegraphics[width=0.5\textwidth]{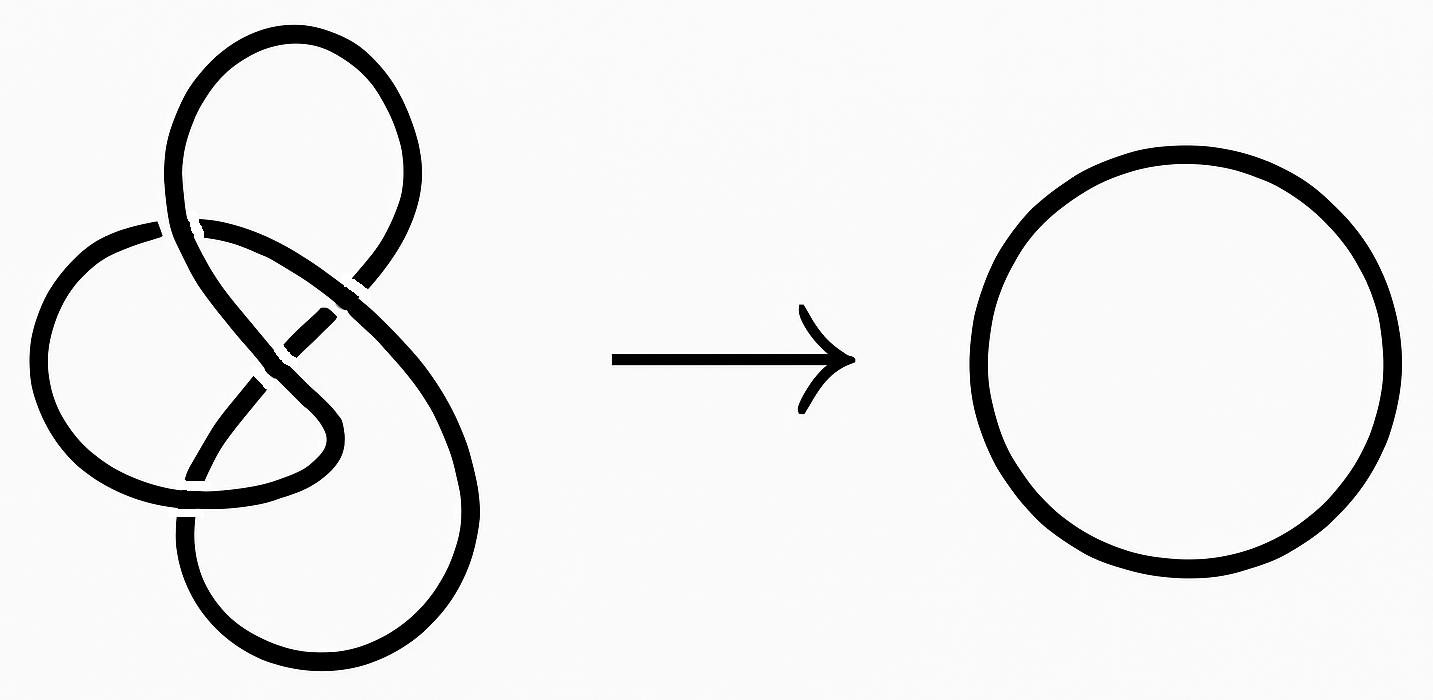}
    \caption{An illustration of ambient isotopy.}
    \label{fig:iso}
\end{figure}

In this context, ``smoothly'' means that the deformation happens gradually and continuously in three-dimensional space. Imagine reshaping a flexible rope without tearing, gluing, or letting any part pass through itself. Mathematically, this is known as a continuous transformation, or \textit{homeomorphism}, of space that preserves the knot’s topological structure. This idea captures the topological essence of knots: we care about the shape of the loop up to continuous deformation, not the exact geometric form. In practice, this means that knots are studied through their behaviour under such flexible transformations rather than by measuring distances or angles. For an illustration of this concept see Figure~\ref{fig:iso}.

\subsection*{Knot Diagrams and Projections}

Because it's hard to study knots directly in 3D, we often look at their 2D \textit{projections}, like casting a shadow of the knot onto a piece of paper. When we do this, we create a \textbf{knot diagram}, which includes information about which strands cross over and under at each crossing.

\begin{tcolorbox}[title=Definition of Knot Diagrams]
A \textbf{knot diagram} is a two-dimensional drawing of a knot that shows all the crossings, with breaks in the lines to indicate which strand goes over or under.
\end{tcolorbox}

\begin{figure}[ht]
\centering
\includegraphics[width=0.5\textwidth]{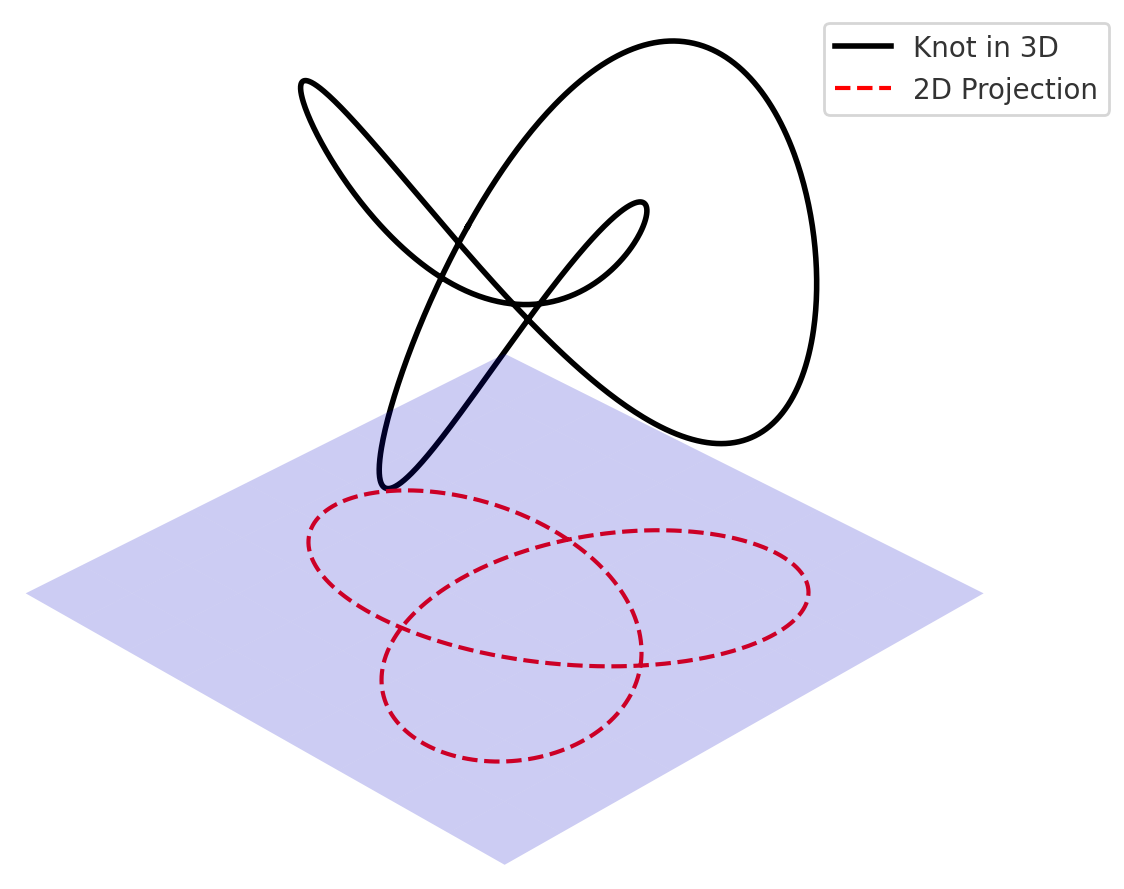}
\caption{A three-dimensional knot (black) is projected onto a two-dimensional plane (light blue), resulting in a knot diagram (red dashed curve). While the 2D projection simplifies the knot's representation, over/under crossing data must be carefully preserved to accurately reflect the knot's structure.}
\label{fig:knot-projection-plane}
\end{figure}

To better understand how knot diagrams are created, we can visualize the process of projecting a three-dimensional knot onto a two-dimensional plane. This projection acts like casting a shadow of the knot onto a plane. However, unlike ordinary shadows, knot diagrams preserve crucial over/under crossing information to maintain the knot’s topological structure. Figure~\ref{fig:knot-projection-plane} illustrates this idea: the black curve represents a 3D knot, the light blue surface is the projection plane, and the red dashed curve shows the resulting two-dimensional projection.

Note that not all planar projections of knots yield valid diagrams. In mathematical knot theory, we require knot diagrams to follow certain conventions that avoid ambiguity:
\begin{itemize}
    \item No triple points: At most two strands may cross at any given location.
\item No tangencies: Strands must cross transversely; they cannot just touch or run alongside each other.
\item No self-intersections except at crossings: The only permitted overlaps are at designated crossings, with clearly marked over/under information.
\end{itemize}

These constraints ensure that diagrams remain interpretable and can be analysed through various topological tools. Any projection violating these conditions must be adjusted by perturbing the knot slightly into a regular diagram.
These diagrams are central to knot theory, as they give us a practical way to study knots, compare them, and apply mathematical tools.

\subsection*{Reidemeister Moves}

\begin{figure}[ht]
\centering
\includegraphics[width=0.5\textwidth]{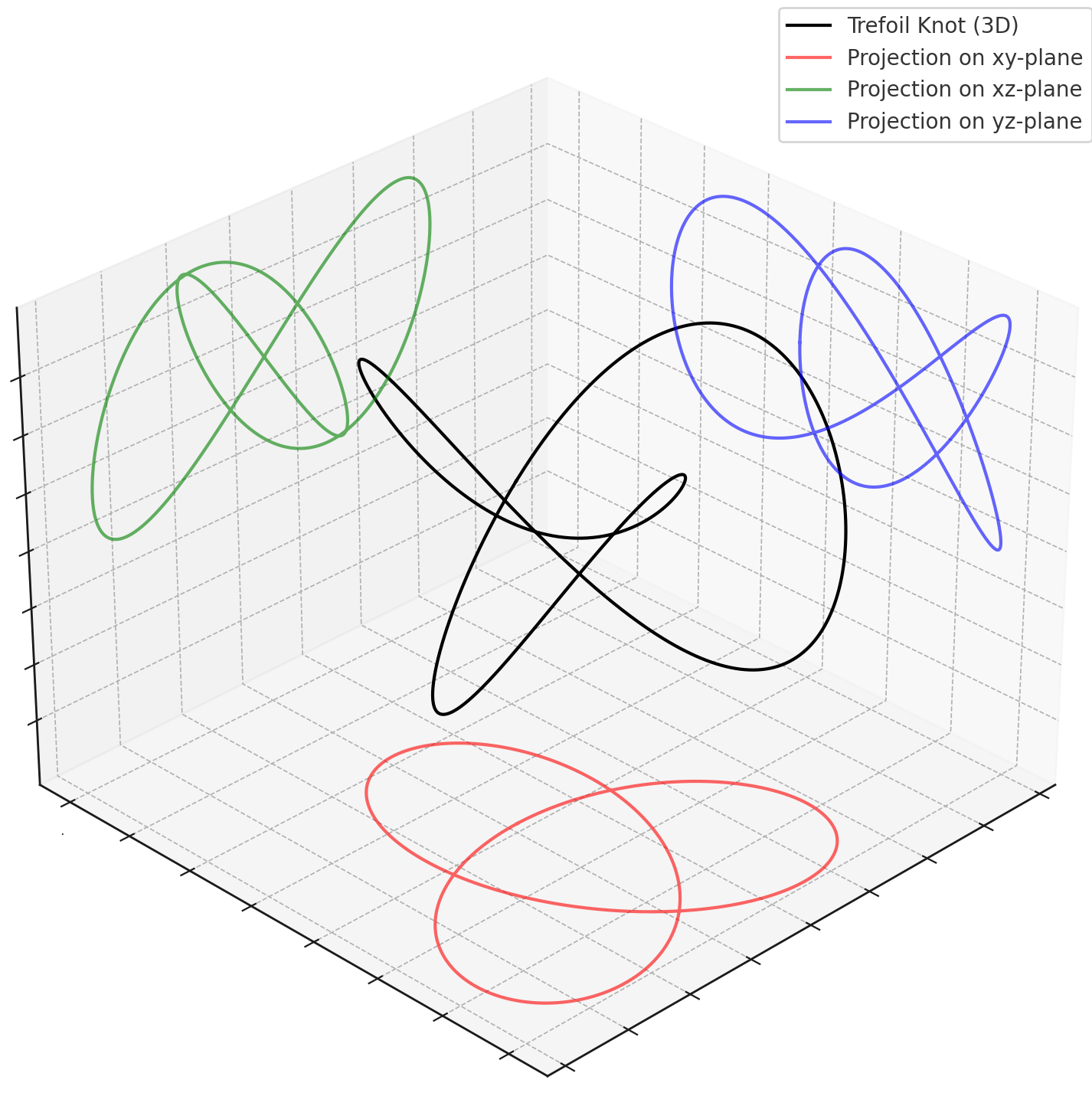}
\caption{The trefoil knot in three-dimensional space (black curve) and its orthogonal projections onto the coordinate planes. The red curve shows the projection on the \(xy\)-plane, green on the \(xz\)-plane, and blue on the \(yz\)-plane. These projections help transition from spatial intuition to two-dimensional knot diagrams.}
\label{fig:trefoil-projections}
\end{figure}

The diagram of a knot is not unique as it depends on the angle or plane from which the knot is projected. A trefoil knot, for instance, may look quite different when projected onto the xy, xz, or yz-plane, as shown in Figure~\ref{fig:trefoil-projections}. Yet all these projections represent the same knot in space. This raises an important question: how can we determine whether two seemingly different diagrams actually correspond to the same underlying knot? To answer this, mathematicians developed a set of simple diagrammatic transformations, known as the {\it Reidemeister moves} \cite{Reid}, that capture precisely when two knot diagrams are topologically equivalent. These moves form the foundation of diagrammatic knot theory and allow us to manipulate, compare, and classify knots using only their two-dimensional representations.

\begin{figure}[ht]
\centering
\includegraphics[width=0.8\textwidth]{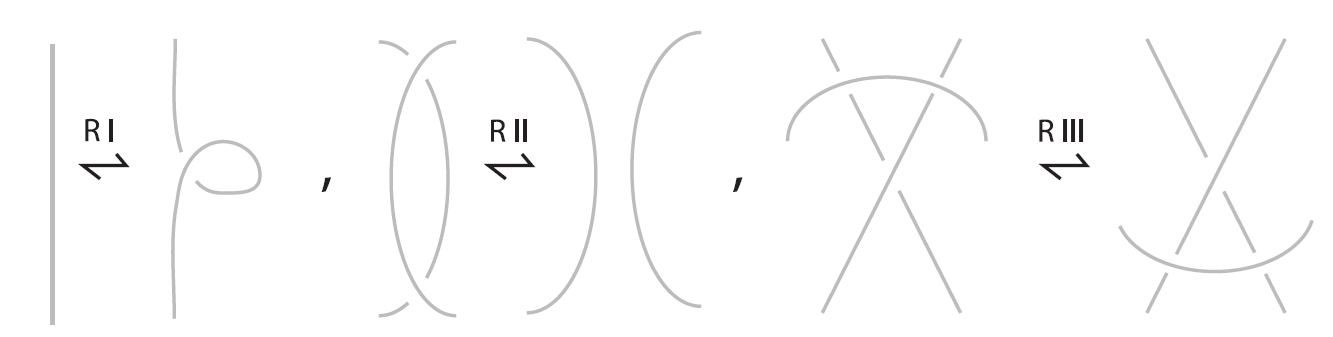}
\caption{The three Reidemeister moves.}
\label{fig:Rmoves}
\end{figure}

\begin{tcolorbox}[title=The Reidemeister Theorem]
There are three Reidemeister moves, each changing a small part of a knot diagram (for an illustration see Figure~\ref{fig:Rmoves}):
\begin{itemize}
\item[] \textbf{RI Move}: Twisting or untwisting a loop.
\item[] \textbf{RII Move}: Sliding one loop over another.
\item[] \textbf{RIII Move}: Moving a strand over or under a crossing.
\end{itemize}
Two knot diagrams represent the same knot if and only if they can be transformed into one another using a sequence of Reidemeister moves.
\end{tcolorbox}

Understanding these moves allows students to explore the idea of knot equivalence directly, using drawings or string models. They offer a beautiful example of how a few simple rules can generate deep mathematical structure.

\medskip

\subsection*{Knot Invariants}

As students work with knots and links, a natural question arises:

\begin{center}
\begin{framed}
{\it How do we know when two knots are really different?}    
\end{framed}
\end{center}

Reidemeister moves help us recognize when two diagrams represent the same knot, but distinguishing different knots requires more powerful tools. This is where \textbf{knot invariants} come into play.

\begin{tcolorbox}[title=Knot Invariants]
A \textbf{knot invariant} is a quantity or property that remains unchanged under ambient isotopy. If two knots have different values for an invariant, they must be different knots.
\end{tcolorbox}

Some basic invariants are intuitive. For example, the \textit{crossing number} of a knot is the smallest number of crossings in any diagram representing that knot. The trefoil knot, for example, cannot be represented with fewer than three crossings, while the Hopf link, a link of two interlocked unknots, has a minimal diagram with two crossings. Since these crossing numbers differ, the trefoil knot is not equivalent to the Hopf link; no sequence of Reidemeister moves can transform one into the other. While the crossing number is a useful starting point, other invariants delve deeper into the structure and classification of knots and links, offering more robust tools for distinguishing between topological types.

These ideas are rich with potential for classroom exploration. They show that behind the intuitive notion of a knot lies a sophisticated mathematical landscape.

\begin{tcolorbox}[title=Big Questions in Knot Theory, colback=blue!5, colframe=blue!50!black]
\begin{itemize}
    \item How many different knots are there?
    \item How can we tell when two knots are really the same?
    \item Are there knots that cannot be untied?
    \item What happens when we tangle more than one loop?
    \item Can we measure how ``knotted'' a knot is?
\end{itemize}
\end{tcolorbox}

These questions can guide students’ curiosity and spark rich discussions. They also help position knot theory as a gateway to deeper mathematical thinking, where visual intuition meets formal reasoning. The next section grounds these questions in educational theory, demonstrating why knot theory is so well-suited for student-centered, exploratory learning.

\section{Theoretical Framework: Constructivism and Visual Mathematics}

Constructivist learning theory emphasizes that knowledge is not passively absorbed but actively constructed by learners through direct experience, experimentation, and reflection  \cite{Piaget1952, Vygotsky1978, Bruner1966}. This approach positions students as agents of their own learning, capable of forming meaningful conceptual understandings when given opportunities to explore, question, and engage with the material. Rather than being told mathematical facts, students encounter problems that challenge them to develop insights, articulate patterns, and test their understanding in a dynamic learning environment.

Knot theory is uniquely suited to this mode of learning. It allows students to manipulate real-world materials, such as ropes or strings, and observe how transformations affect the structure of a knot. Through trial and error, they discover that not all transformations are equal, leading to the development of rules and invariants. The realization that certain moves preserve a knot while others fundamentally alter it introduces the concept of mathematical structure in an intuitive and memorable way. 

The tactile and visual nature of knot theory also makes it ideal for students who may struggle with more symbolic or algebraic topics in traditional mathematics. In knot theory, the abstract is rendered concrete: knots can be drawn, touched, and transformed. For example, when students investigate knot equivalence through the lens of Reidemeister moves, they engage directly with mathematical reasoning by manipulating diagrams and making predictions about their outcomes.

Visual-spatial reasoning plays a central role. Students must interpret diagrams, imagine movements in three-dimensional space, and translate between physical knots and their two-dimensional representations \cite{Arcavi2003, Duval1999}. With the addition of visual figures in this paper (illustrating projections, isotopies, and diagrammatic transformations), students can better grasp the inherently geometric nature of knot theory. These tasks promote cognitive flexibility, spatial awareness, and a deeper understanding of how abstract concepts can be grounded in intuition.

Furthermore, the collaborative nature of many knot theory activities, such as classifying knots, discovering equivalences, or creating artistic representations, aligns well with the social aspects of constructivist pedagogy \cite{Vygotsky1978, Bruner1966}. Learning in groups encourages dialogue, explanation, and refinement of ideas, making the mathematical experience richer and more memorable. Peer discussion helps students articulate their thinking, confront misconceptions, and co-construct knowledge.

\smallbreak

In the context of a knot theory lesson, constructivist principles are realized through:
\begin{itemize}
\item \textbf{Concrete manipulations}: tying and untying physical knots, using colored string to highlight crossings, and interacting with virtual knot simulators.
\item \textbf{Visual-spatial reasoning}: interpreting and creating 2D projections of 3D objects, recognizing equivalent diagrams, and visualizing Reidemeister moves.
\item \textbf{Collaborative problem-solving}: working in teams to classify knots, challenge each other with equivalence puzzles, or recreate complex knots from minimal descriptions.
\item \textbf{Intuitive abstraction}: formulating rules based on empirical observations, such as developing informal criteria for when two knots are the same or different.
\end{itemize}

Ultimately, this approach cultivates skills central to mathematical maturity: inquiry, persistence in problem-solving, creative representation of abstract ideas, and effective communication of reasoning \cite{Bruner1966, Arcavi2003}. By learning through exploration, students not only grasp mathematical ideas more deeply, they also come to see themselves as capable of doing mathematics in meaningful, authentic ways \cite{devlin, Hmelo2004}. Knot theory becomes not just a topic to be studied, but a medium through which students develop a mathematical identity \cite{Presmeg2006}.

\section{The VIBE Framework}

Designing engaging and conceptually rich classroom experiences around knot theory requires more than just a collection of activities. It demands a coherent pedagogical strategy. To this end, we introduce the VIBE framework, a structured yet adaptable approach for translating the constructivist and visual foundations discussed in Section 2 into practice.

The name VIBE (Visual, Inquiry-based, Braided, Embedded) captures the rhythm and flow of effective mathematical engagement, especially within the rich and tactile domain of knot theory. Each pillar represents a key dimension of meaningful learning, but it is their interplay that brings the framework to life. VIBE is not meant to be a rigid sequence of steps; instead, it is a flexible model that supports active, creative, and collaborative exploration of mathematical ideas.

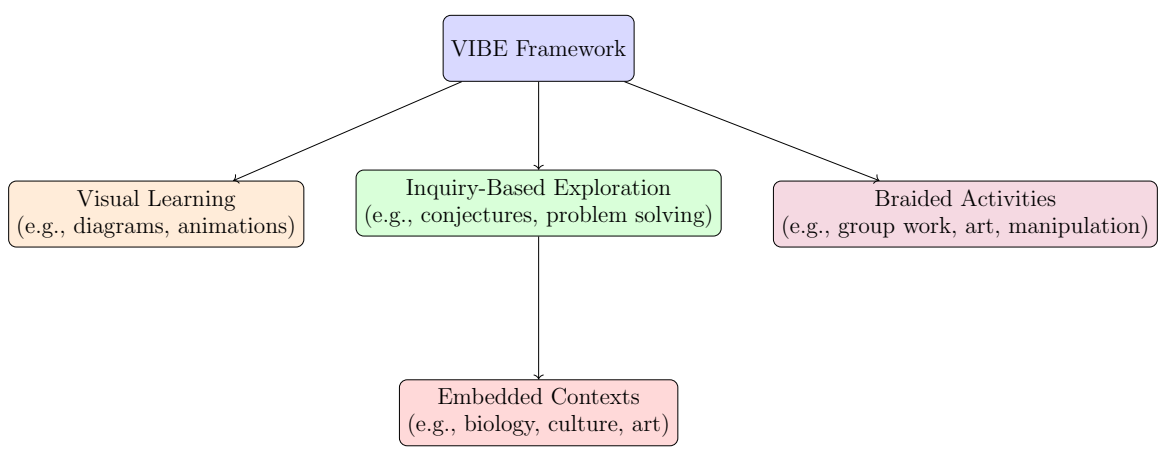
\begin{figure}[ht]
\centering
\resizebox{0.95\textwidth}{!}{%
\begin{tikzpicture}[node distance=1.6cm, every node/.style={rectangle, draw, rounded corners, align=center, minimum height=1.2cm}]
\node (vibe) [fill=blue!15] {VIBE Framework};

\node (visual) [below left=1.8cm and 2.5cm of vibe, fill=orange!15] {Visual Learning\\(e.g., diagrams, animations)};
\node (inquiry) [below=of vibe, fill=green!15] {Inquiry-Based Exploration\\(e.g., conjectures, problem solving)};
\node (braided) [below right=1.8cm and 2.5cm of vibe, fill=purple!15] {Braided Activities\\(e.g., group work, art, manipulation)};
\node (embedded) [below=2.6cm of inquiry, fill=red!15] {Embedded Contexts\\(e.g., biology, culture, art)};

\draw[->] (vibe) -- (visual);
\draw[->] (vibe) -- (inquiry);
\draw[->] (vibe) -- (braided);
\draw[->] (inquiry) -- (embedded);

\end{tikzpicture}%
}
\caption{Schematic overview of the VIBE pedagogical framework.}
\label{fig:vibe-chart}
\end{figure}

The components of the VIBE framework are presented below.

\begin{itemize}
    \item \textbf{V -- Visual Learning}: Visual representations play a foundational role in students’ mathematical thinking and are especially effective in topological contexts \cite{Presmeg2006, Arcavi2003}. Knot theory is inherently visual. From the moment students hold a piece of string to when they begin drawing knot diagrams or exploring digital models, they engage with mathematics in a spatial and geometric way. Visual representations help bridge the gap between concrete experiences and abstract thinking. Teachers can use animations to demonstrate knot transformations, color-coded diagrams to introduce invariants, and software like KnotPlot to allow students to manipulate knots in 3D.
    \item \textbf{I -- Inquiry-Based Exploration}: The inquiry-based pillar of VIBE is strongly influenced by the success of problem-based learning models \cite{hmelo2000, Hmelo2004, Bell2010}. Rather than presenting rules and definitions from the outset, this pillar encourages students to ask their own questions and investigate patterns. For example, students might wonder whether two knots are the same or what operations preserve a knot’s structure. Teachers act as facilitators, prompting students to form conjectures, test their ideas, and build a deeper understanding through discovery. This process mirrors authentic mathematical research and helps students develop critical thinking skills.

\item \textbf{B -- Braided Activities}: The ``braid'' metaphor captures the interactive, intertwined nature of effective classroom activities. In this pillar, students engage in hands-on tasks that weave together different aspects of knot theory. These may include drawing knots, simplifying diagrams using Reidemeister moves, classifying knots using colorability rules, or creating their own artistic or functional knots. Group activities and collaborative challenges foster communication, teamwork, and shared insight.

\item \textbf{E -- Embedded Contexts}: Mathematical ideas become more meaningful when they are situated in rich, interdisciplinary contexts. Knot theory offers numerous natural connections to other subjects: DNA replication in biology, network routing in computer science, weaving and design in visual arts, and even historical uses of knots in navigation and craftsmanship. By embedding knot theory in these contexts, teachers make abstract mathematics relevant and inspiring. Incorporating cultural representations of knots supports broader global citizenship goals in mathematics education \cite{UNESCO2015, Nussbaum2006}.

\end{itemize}

\smallbreak

Together, the VIBE pillars form a versatile instructional model that empowers both teachers and students. Educators can adapt them to fit different classroom settings, learning goals, and time constraints. The framework supports diverse learners by balancing tactile, visual, and conceptual engagement, and it invites creativity through open-ended exploration.

What distinguishes VIBE from traditional lesson planning models is its grounding in authentic mathematical inquiry. It mirrors the practices of mathematicians: visualizing, asking questions, collaborating, and connecting ideas across domains. Moreover, because it is rooted in familiar experiences (drawing, tying knots, storytelling, and problem-solving) VIBE lowers the threshold for participation while raising the ceiling for mathematical depth.

VIBE also serves as a bridge between curriculum and context. Whether implemented as a short enrichment unit, a math club theme, or part of a cross-disciplinary project, the framework promotes an active, connected, and student-centered approach to learning mathematics through knots. In doing so, it helps cultivate mathematical maturity and curiosity in ways that resonate beyond the classroom.

\begin{remark}\rm
While the acronym VIBE has been used in educational literature related to equity in higher education, our usage, Visual‑Inquiry‑Braided‑Embedded, is distinct and original to this pedagogical application in knot theory.
\end{remark} 

Educators can adapt the VIBE pillars to fit different classroom settings, student backgrounds, and curricular goals. Whether implemented as a short enrichment unit, a math club theme, or part of a cross-disciplinary project, the VIBE framework promotes an active, connected, and student-centered approach to learning mathematics through knots.

\subsection*{Visual Summary of Pedagogical Dimensions}

\begin{figure}[ht]
\centering
\includegraphics[width=0.85\textwidth]{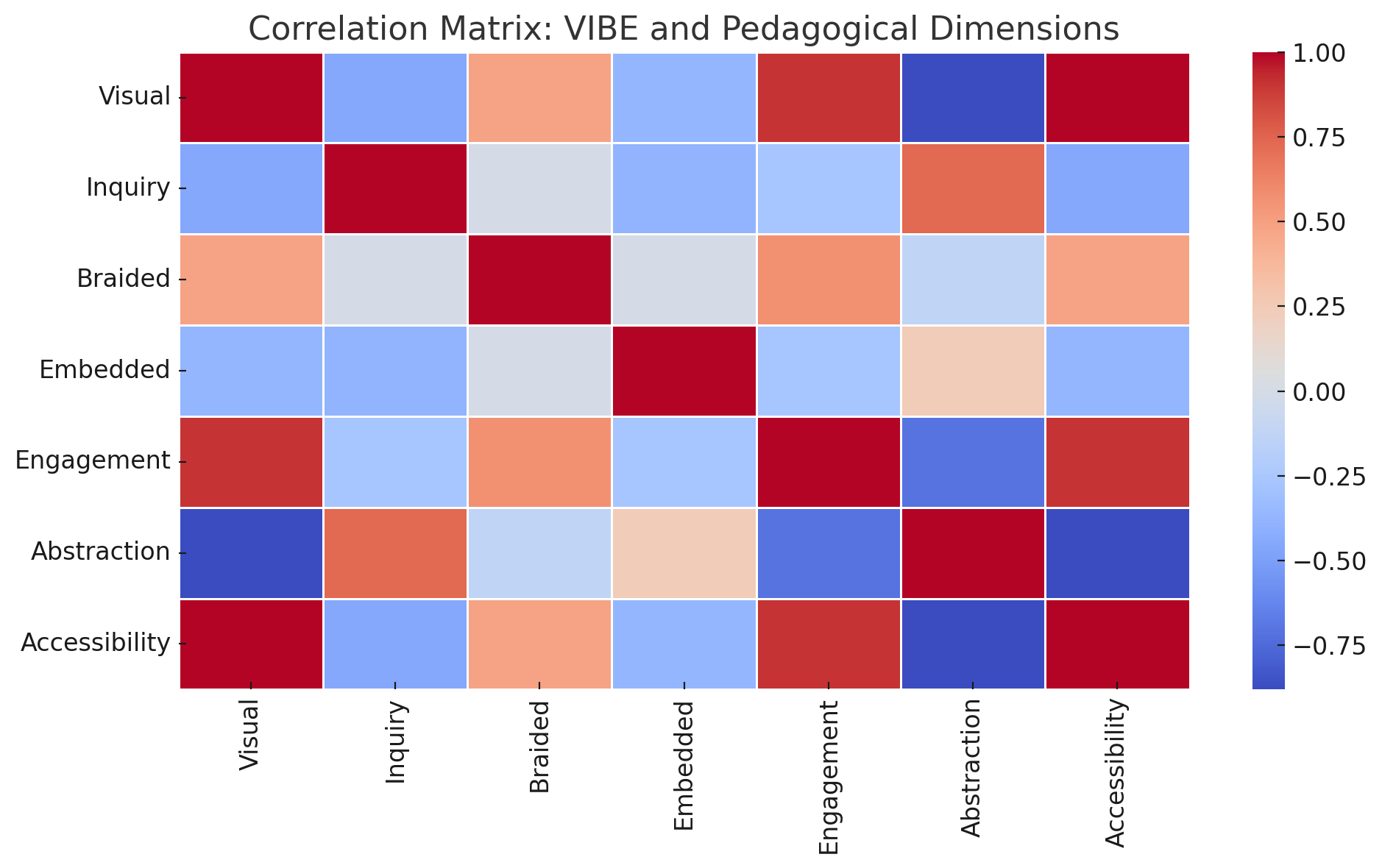}
\caption{Illustrative heatmap of VIBE pillars and their qualitative alignment with key pedagogical outcomes. Darker shades indicate stronger associations.}
\label{fig:heatmap}
\end{figure}

To visually synthesize the VIBE framework’s pedagogical reach, we present a qualitative heatmap aligning its four pillars with central educational goals: engagement, abstraction, accessibility, and collaboration. While no numeric values are assigned, the colour intensity represents the estimated strength of each association based on theoretical reasoning and classroom practice (for an illustration see Figure~\ref{fig:heatmap}).

For instance, visual learning is intuitively linked to both accessibility and engagement, as it lowers entry barriers and sparks curiosity. Inquiry-based learning fosters abstraction by encouraging students to formulate and test hypotheses. Braided activities emphasize collaboration through interactive, hands-on tasks, and embedded contexts enhance engagement by making mathematics meaningful across disciplines and cultures.

\begin{tcolorbox}[colback=green!5!white, colframe=green!50!black, title=Interpretation and Future Use]
This heatmap is not based on empirical data but serves as a design tool for educators and researchers. It invites further validation through classroom experimentation, student feedback, or observational studies aimed at measuring the actual impact of each pillar in different learning environments.
\end{tcolorbox}

\begin{figure}[ht]
\centering
\includegraphics[width=0.6\textwidth]{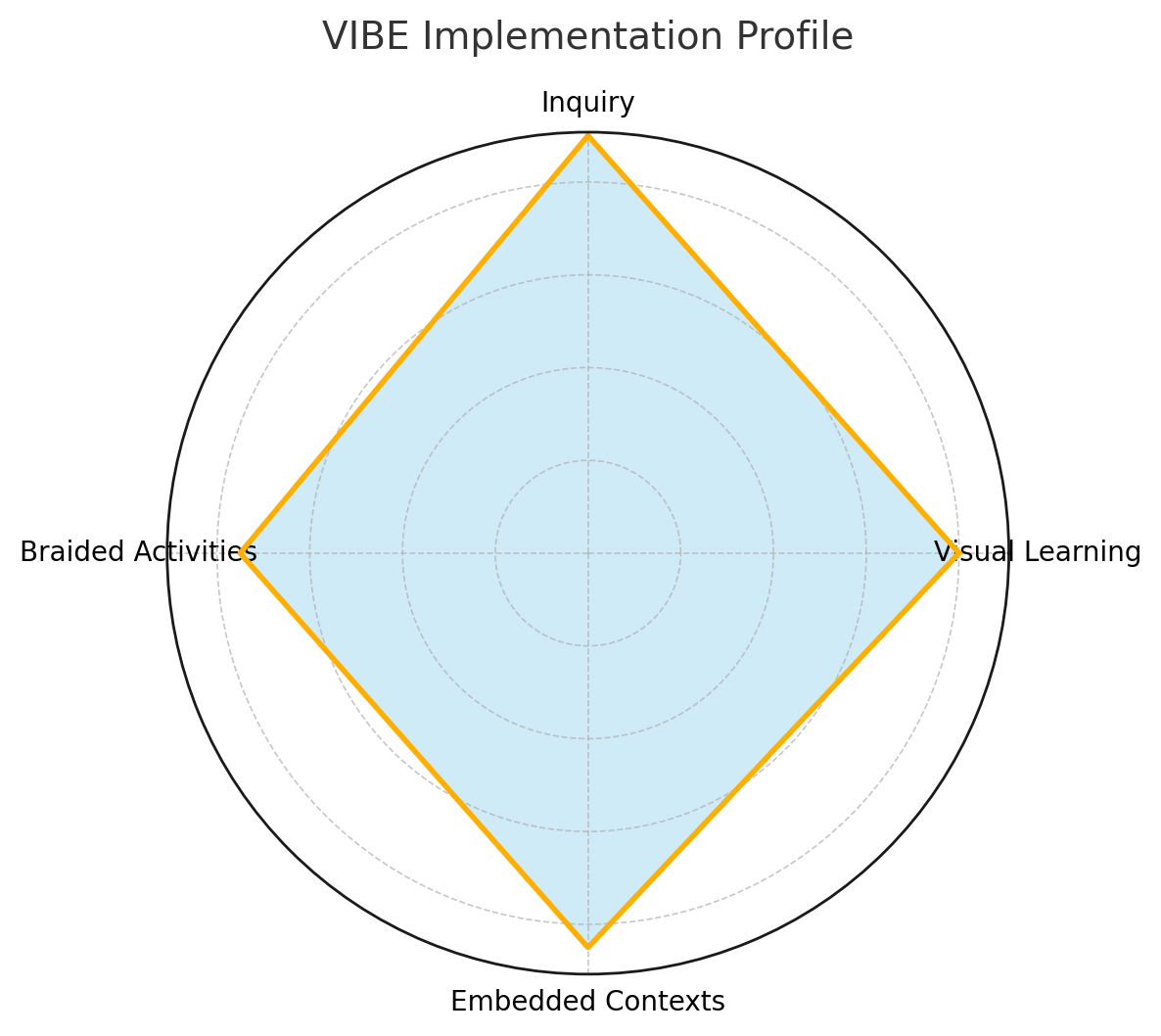}
\caption{The relative strength of each VIBE pillar in a hypothetical classroom or curriculum.}
\label{fig:radar}
\end{figure}

To complement the heatmap, we include a radar chart that visually summarizes the implementation balance across the four pillars of the VIBE framework (see Figure~\ref{fig:radar}). While all pillars are well-represented, this visualization helps identify areas for future emphasis or curricular strengthening. Such profiles can be adapted by educators to reflect their own classroom contexts or institutional priorities.

\subsection*{Strategic Evaluation of the VIBE Framework}

In order to assess the broader pedagogical potential of the VIBE framework beyond individual classroom activities, we present a strategic SWOT (Strengths, Weaknesses, Opportunities, Threats) analysis. This evaluation focuses on VIBE as an instructional design model, considering its alignment with contemporary educational paradigms, its potential for interdisciplinary integration, and its theoretical coherence. The analysis complements the empirical insights provided by the radar chart and heatmap, offering a reflective lens on the model’s adaptability and strategic viability.

\begin{tcolorbox}[colback=gray!5!white, colframe=black!75!white, title=SWOT Analysis of the VIBE Framework (Pedagogical Model)]
\begin{tabularx}{\textwidth}{X|X}
\textbf{Strengths} &
\textbf{Weaknesses} \\
\begin{itemize}[leftmargin=*, label=--]
    \item Visually engaging and spatially intuitive
    \item Supports constructivist and inquiry-based learning
    \item Promotes inclusivity and differentiated instruction
    \item Grounded in authentic mathematical reasoning
\end{itemize}
&
\begin{itemize}[leftmargin=*, label=--]
    \item May require significant teacher training
    \item Risk of being reduced to superficial visual learning
    \item Less aligned with standardized assessment cultures
    \item Limited existing curricular support
\end{itemize}
\\
\hline
\textbf{Opportunities} &
\textbf{Threats} \\
\begin{itemize}[leftmargin=*, label=--]
    \item Adaptation to other STEM and non-STEM subjects
    \item Integration with digital and VR technologies
    \item Application in project-based and cross-curricular contexts
    \item Publication of modular resources and toolkits
\end{itemize}
&
\begin{itemize}[leftmargin=*, label=--]
    \item Resistance from traditionalists in math education
    \item Institutional barriers to curricular innovation
    \item Overshadowed by emerging AI-focused curricula
    \item Risk of inequitable access across school systems
\end{itemize}
\end{tabularx}
\end{tcolorbox}

\noindent This analysis reinforces the pedagogical promise of the VIBE framework while acknowledging areas where support, development, or advocacy may be needed for successful implementation. It provides a strategic foundation for future research and adaptation efforts across disciplines.

\section{Sample Activities and Lessons}

The activities in this section are designed to reflect the VIBE framework in practice, offering a structured yet flexible approach for incorporating knot theory into the classroom. The lessons are adaptable for different age groups, time frames, and educational contexts. Each day’s plan includes concrete tasks, collaborative components, and opportunities for inquiry and creative expression.

\subsection*{Mini-Unit Plan (3--5 Lessons)}

\subsubsection*{Day 1: What is a Knot?}\,

\smallbreak

\noindent \textbf{Objective}: Introduce students to the concept of mathematical knots and develop intuition through hands-on experience.
\smallbreak

\noindent \textbf{Activities}:
\begin{itemize}
    \item Begin with a discussion: ``What is a knot?'' Have students list examples from daily life (shoelaces, climbing rope, sailing knots, etc.).
    \item Distribute ropes or strings and ask students to tie various knots. Can they recreate what they use to tie their shoes? Can they invent a new knot?
    \item Introduce the definition of a mathematical knot: a closed loop (no ends) in three-dimensional space.
    \item Students attempt to draw the knot they tied as a 2D diagram.
    \item Challenge: Give students pairs of knot diagrams (e.g., unknot vs. trefoil) and ask them to determine which are the same and which are different.
\end{itemize}

\subsubsection*{Day 2: Reidemeister Moves and Diagrammatic Equivalence}\,

\smallbreak
\noindent \textbf{Objective}: Explore the equivalence of knots through allowable transformations.

\smallbreak
\noindent \textbf{Activities}:
\begin{itemize}
    \item Introduce the three Reidemeister moves with physical demonstrations (e.g., loops and twists in string) and diagrammatic versions.
    \item Provide students with pairs of diagrams and ask them to try converting one into the other using a sequence of Reidemeister moves.
    \item Group activity: ``Knot Surgery'': each team receives a complex diagram and must simplify or manipulate it using allowed moves to reach a target diagram.
    \item Discuss the idea of knot equivalence and the challenge of determining when two knots are the same.
\end{itemize}

\subsubsection*{Day 3: Tricolorability and Knot Invariants}\,

\smallbreak
\noindent \textbf{Objective}: Discover a simple knot invariant and use it to classify knots.

\smallbreak
\noindent \textbf{Background: What is Tricolorability?}
Tricolorability is a simple way to distinguish certain knots using colouring rules \cite{adams2004}. Given a knot diagram, we try to colour each arc (a segment between crossings) using exactly three colours, following a specific rule: at every crossing, the three arcs must either all have the same colour or all have different colours. If this is possible, the knot is said to be tricolorable. This rule is more than a colouring game, it is a mathematical \textit{invariant}, meaning that it stays the same even if we change the diagram using Reidemeister moves. Therefore, if one knot is tricolorable and another isn’t, they cannot be the same knot.

\smallbreak
\noindent \textbf{Activities}:
\begin{itemize}
    \item Explain tricolorability: students must colour arcs of a knot diagram using three colours with the rule that at each crossing either all three colours appear or only one colour appears.
    \item Apply this to the unknot, trefoil, and figure-eight diagrams. Which are tricolorable?
    \item Discuss why this property is useful: if two diagrams are truly equivalent, they must share this invariant.
    \item Create a classroom chart summarizing which knots are tricolourable.
\end{itemize}

\subsubsection*{Day 4: Knots in the Real World}\,

\smallbreak
\noindent \textbf{Objective}: Connect knot theory to biology, art, and culture.
\smallbreak

\noindent \textbf{Background: Knot Theory in Biology.}
In molecular biology, long DNA strands often become entangled within the limited space of a cell nucleus. Enzymes such as topoisomerases and recombinases modify these structures by cutting and rejoining strands, processes that can be modelled using knot theory \cite{flapan}. 

\begin{figure}[ht]
\centering
\includegraphics[width=0.4\textwidth]{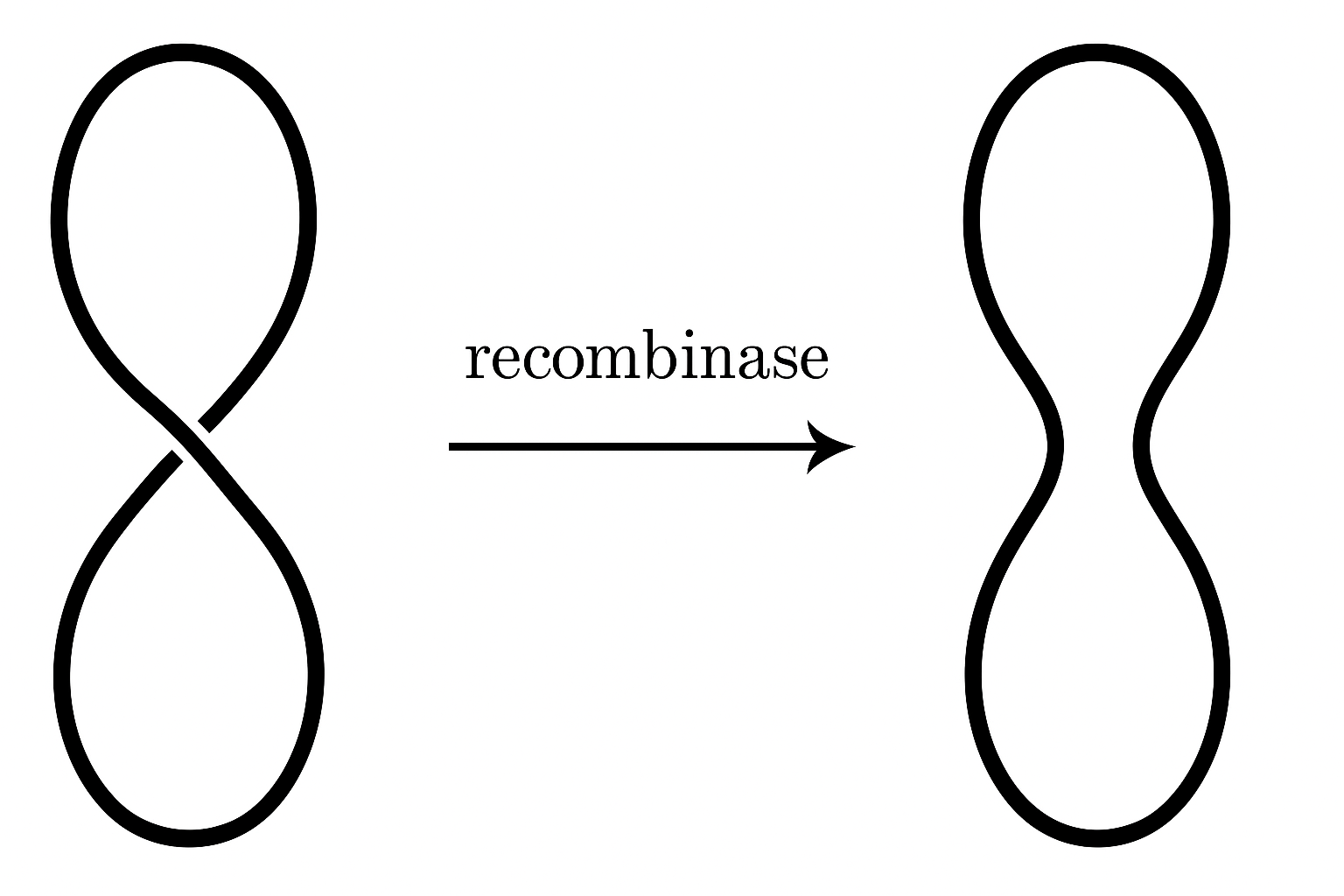}
\caption{A schematic of site-specific recombination: a segment of DNA initially crossing over another is modified by a recombinase enzyme. This is over-simplified and meant as an educational illustration only.}
\label{fig:recombinase}
\end{figure}

As shown in Figure~\ref{fig:recombinase}, the action of a recombinase can change the crossing information in a DNA strand, effectively transforming one knot or link type into another. By representing DNA as a closed curve (or a tangle of multiple curves) biologists can use knot invariants to distinguish molecular configurations, analyse folding processes, and predict the effects of enzymatic transformations on genome structure. This interdisciplinary bridge offers a compelling example for students: abstract mathematics becomes a lens through which concrete phenomena in modern biology can be understood and explored.

\smallbreak
\noindent \textbf{Activities}:
\begin{itemize}
    \item Presentation or video: DNA supercoiling and knotted proteins; how knot theory helps biologists understand molecular behavior.
    \item Art activity: Design a knot-based pattern for jewelry, textiles, or logos. Provide examples from Celtic knots or modern graphic design.
    \item Short reflective writing: ``Where do you think knot theory could be applied in the future?''
\end{itemize}

\subsubsection*{Day 5 (Optional): Final Project and Gallery Walk}\,

\smallbreak
\noindent \textbf{Objective}: Consolidate learning through presentation and peer assessment.

\smallbreak
\noindent \textbf{Activities}:
\begin{itemize}
    \item Students prepare posters or digital slides summarizing what they learned about a particular knot (e.g., trefoil, figure-eight, or their own creation).
    \item Include elements such as history, classification, invariants, visual diagrams, and potential applications.
    \item Conduct a gallery walk: students display their work and peer-review each other’s projects using comment cards.
\end{itemize}

\subsection*{Supplementary Activities}

\begin{description}
    \item[Activity 1: Knot Gallery Walk] Post diagrams of different knots around the classroom, each labeled with a number. Students classify each knot and justify their reasoning. A class-wide discussion follows.

    \item[Activity 2: Create Your Own Knot Invariant] Students invent and test their own methods of distinguishing knots. Possibilities include counting crossings, noting handedness, or using color-based systems. Groups test the robustness of their invariants.

    \item[Activity 3: Knot and Braid Theater] A kinesthetic activity where students act as strands of a braid or knot, moving with ropes to simulate over/under crossings. Helps internalize the dynamics of topology.

    \item[Activity 4: Build a Knot Museum] Students create a ``Knot Museum'' corner with posters, physical models, and artistic or historical exhibits. Other classes or parents are invited to tour the exhibit.
\end{description}

Each of these activities exemplifies the VIBE pillars, promoting visual learning, inquiry, collaboration, and interdisciplinary relevance. Differentiation is built in: advanced learners can explore algebraic invariants (e.g., the Jones polynomial), while others focus on visualization and classification.

\begin{remark}\rm
It is worth mentioning that game-based approaches have gained traction in knot pedagogy. For example, Henrich et al. \cite{henrich2010} introduced interactive knot-shadow games such as ``To Knot or Not To Knot'', which use playful challenges to deepen spatial intuition and foster strategic reasoning. More recently, Henrich and collaborators developed knot-themed games such as KNOTRIS and the Arc Unknotting Game, blending problem-solving with exploratory play \cite{henrichknotris, henricharcunknot}.
\end{remark}

\smallbreak
\subsection*{Evaluation of the Pedagogical Design}

To further evaluate the pedagogical design of the proposed knot theory activities, we now present a heatmap that synthesizes two key educational dimensions: cognitive demand and accessibility. Each activity is qualitatively placed in relation to these dimensions using a colour-coded scheme, where green indicates high accessibility, yellow represents medium levels, and red reflects greater complexity or lower accessibility. This visualization serves multiple purposes: it helps educators select or adapt activities based on their classroom context, it reveals the balance between challenge and inclusiveness, and it underscores how VIBE-aligned tasks can be differentiated to meet diverse learner needs.

Activities marked with deeper red shades under ``Cognitive Demand'' signal rich, high-level reasoning opportunities (e.g., creating knot invariants), while those shaded green under ``Accessibility'' highlight low-barrier entry points ideal for engaging a broad range of students (e.g., artistic knot design or braid theater). By visually mapping these dimensions, the heatmap allows educators to scaffold instruction, mix and match tasks across levels, and ensure all learners find pathways into mathematical exploration.

\definecolor{Low}{RGB}{198,239,206}
\definecolor{Med}{RGB}{255,235,156}
\definecolor{High}{RGB}{255,199,206}

\begin{table}[h!]
\centering
\small
\renewcommand{\arraystretch}{1.3}
\begin{tabular}{|l|c|c|}
\hline
\rowcolor{gray!20}
\textbf{Activity} & \textbf{Cognitive Demand} & \textbf{Accessibility} \\
\hline
Drawing knots & \cellcolor{Med} Medium & \cellcolor{Low} High \\
Reidemeister move puzzles & \cellcolor{High} High & \cellcolor{Med} Medium \\
Tricoloring & \cellcolor{Med} Medium–High & \cellcolor{Med} Medium \\
Artistic knot design & \cellcolor{Low} Low & \cellcolor{Low} High \\
Knot classification tasks & \cellcolor{High} High & \cellcolor{Med} Medium \\
DNA connections & \cellcolor{Med} Medium & \cellcolor{Low} High \\
Braid theater & \cellcolor{Med} Medium & \cellcolor{Low} High \\
Create-your-own invariant & \cellcolor{High} Very High & \cellcolor{High} Low–Medium \\
\hline
\end{tabular}
\caption{Activity Heatmap: Balancing Cognitive Demand and Accessibility}
\end{table}

\begin{figure}[ht]
\centering
\includegraphics[width=0.98\textwidth]{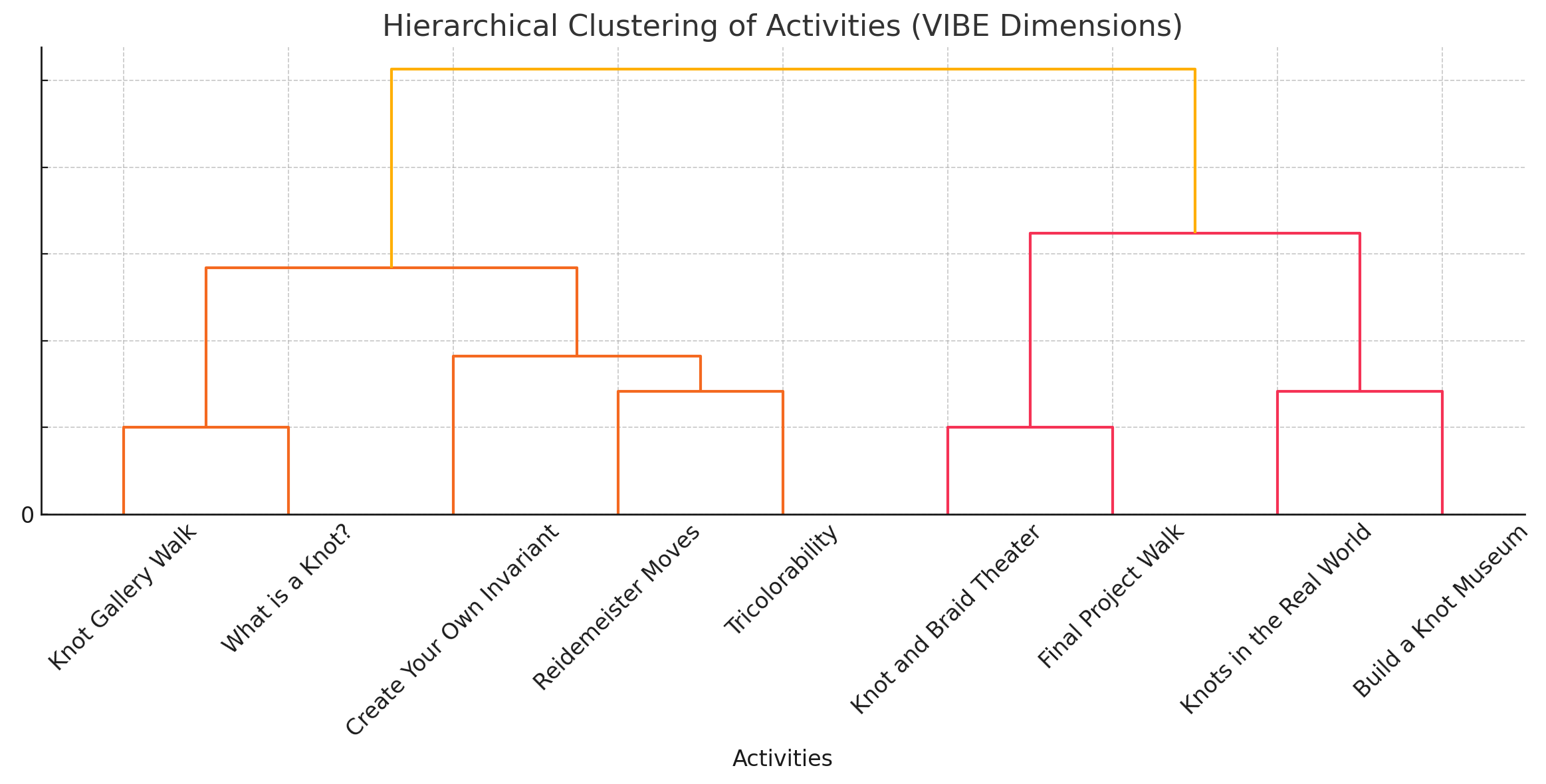}
\caption{Hierarchical clustering of sample activities based on their VIBE profiles. Each activity is scored on the four pedagogical dimensions—Visual, Inquiry, Braided (collaborative), and Embedded (contextual). The dendrogram illustrates clusters of pedagogical similarity, offering insight into how lessons can be grouped or sequenced effectively.}
\label{fig:cluster}
\end{figure}

To illustrate the pedagogical relationships among the proposed knot theory activities, we also present a {\it conceptual dendrogram} based on the VIBE dimensions (see Figure~\ref{fig:cluster}). Rather than being derived from empirical data, this visual grouping reflects intuitive similarities in structure, objectives, and learning modalities. Activities that emphasize hands-on engagement, visual reasoning, and real-world connections naturally cluster together, while those focused on abstraction or theoretical depth form their own branches. This representation can guide educators in designing coherent lesson sequences and balancing different instructional modes.

\begin{tcolorbox}[colback=green!5!white, colframe=green!50!black, title=Interpretation Note]
This dendrogram is a conceptual visualization. It highlights how certain knot theory activities share similar pedagogical profiles, helping educators to design balanced, interconnected lessons that reinforce multiple learning dimensions.
\end{tcolorbox}

At the leftmost branch of the dendrogram, Knot Gallery Walk and What is a Knot? form the closest cluster. These are highly accessible, introductory activities that emphasize visual engagement and foundational understanding, aligning strongly with the Visual and Embedded dimensions of the VIBE framework.

Just above them, Create Your Own Invariant, Reidemeister Moves, and Tricolorability form a nested cluster, reflecting a shared focus on Inquiry and Abstraction. These tasks challenge students to explore mathematical properties, test equivalences, and engage with deeper structural reasoning, making them suitable for more advanced or reflective stages of the learning sequence.

On the right side of the dendrogram, Knot and Braid Theater, Final Project Walk, Knots in the Real World, and Build a Knot Museum form a second major grouping. These activities score highly on the Braided and Embedded dimensions, as they are collaborative, contextual, and often cross-disciplinary. Their clustering suggests that they collectively support applied, performative, or project-based modes of instruction.

This clustering supports differentiated curriculum design: educators can select activities from each group depending on whether their focus is on conceptual exploration, introductory understanding, or real-world applications.

The dendrogram serves as a map of pedagogical design, highlighting which activities reinforce each other, how they may be sequenced, and where there is room for thematic coherence or differentiated instruction. It can also be used to scaffold units, ensuring a balanced coverage of all four pillars in the VIBE framework.

\section{Implementation and Assessment Strategy}

Bringing knot theory into the secondary classroom requires thoughtful implementation and a flexible, student-centered approach to assessment. While the activities and framework outlined in this paper are designed to be accessible and engaging, their success depends on alignment with broader educational goals and classroom realities.

\subsection*{Implementation Guidelines}

\begin{enumerate}
    \item \textbf{Modular Design}: Teachers are encouraged to adopt the lessons as a stand-alone mini-unit, a supplemental enrichment module, or a cross-disciplinary project. The activities can be scaled for 45-minute periods or extended sessions, depending on time availability.
    
    \item \textbf{No Prerequisite Knowledge}: The material assumes no prior exposure to topology or higher-level mathematics. Teachers can scaffold new concepts through concrete experiences, vocabulary development, and peer discussion.
    
    \item \textbf{Interdisciplinary Opportunities}: Knot theory lends itself to integration with biology (e.g., DNA knotting), art (e.g., knot patterns and visual symmetry), and computer science (e.g., network topology). Collaboration with teachers from other disciplines can enhance engagement and deepen relevance.
    
    \item \textbf{Low-Cost Materials}: Most activities require only simple materials: string, markers, paper, scissors, and access to printed or projected diagrams. Optional software tools (such as KnotPlot \cite{plot}) can enrich the experience but are not required.
    
    \item \textbf{Inclusive and Differentiated Instruction}: Knot theory is accessible to a wide range of learners. Visual learners benefit from diagrams; kinesthetic learners thrive on physical manipulation; verbal learners gain from group explanation. Advanced students can explore algebraic knot invariants or braid group representations.
\end{enumerate}

\subsection*{Assessment Strategy}

Given the exploratory and visual nature of the module, traditional testing is neither required nor optimal. Instead, a variety of formative and summative strategies can be used to evaluate learning:

\begin{itemize}
    \item \textbf{Learning Journals}: Students maintain a reflective notebook documenting their processes, insights, and challenges. Prompts can include: \textit{``Describe how you determined whether two knots were the same''}, or \textit{``Explain your strategy for tricoloring the trefoil''.}
    
    \item \textbf{Peer Teaching and Presentations}: Students explain knot diagrams or invariants to peers, fostering communication and consolidation of understanding. Group presentations on specific knots or applications (e.g., DNA, Celtic art) can serve as assessment artifacts.
    
    \item \textbf{Creative Demonstrations}: Posters, visual essays, or digital media projects that showcase student learning through diagrams, storytelling, or artistic interpretation. Example task: \textit{``Create a poster describing the figure-eight knot, its diagram, whether it is tricolorable, and how it appears in real life''.}
    
    \item \textbf{Gallery Walks and Peer Review}: Students circulate to view each other’s work and provide structured feedback using rubrics or guided questions. Sample peer prompt: \textit{“What did you learn from this knot? What could be added or clarified?”}
    
    \item \textbf{Exit Tickets and Concept Maps}: At the end of each lesson, students respond to quick prompts: \textit{“What is one question you still have?”} or \textit{“What was surprising today?”} Concept maps help visualize how terms like \textit{knot}, \textit{invariant}, \textit{Reidemeister move}, and \textit{equivalence} connect.
\end{itemize}

\subsection*{Rubric Example: Tricolorability Project}

\begin{itemize}
    \item \textbf{Clarity of Explanation} (0--5): Student clearly explains the rules of tricolorability and applies them correctly.
    \item \textbf{Visual Accuracy} (0--5): Diagrams are neat, labeled, and accurately show crossings and colorings.
    \item \textbf{Insight and Reflection} (0--5): Student offers personal insight into the process, challenges faced, or mathematical connections.
    \item \textbf{Creativity or Real-World Connection} (0--5): The project shows creativity or links knot theory to broader topics.
\end{itemize}

\smallbreak

Assessment in this framework is not about right or wrong answers, but about growth in reasoning, visualization, persistence, and the ability to communicate ideas clearly. By aligning with constructivist pedagogy, assessment becomes a tool for reflection and deeper learning rather than a gatekeeping mechanism.

Teachers may optionally collect anonymized student work over time to inform broader educational discussions or to contribute to research on visual mathematics and alternative curriculum development.

\subsection*{Classroom Snapshots: Simulated Dialogues}

\begin{tcolorbox}[title=Classroom Snapshot 1: Discovering the Trefoil, colback=gray!5, colframe=gray!75!black]
\textbf{Teacher}: So, what do you notice about this knot? \
\textbf{Student}: It loops three times, kind of like a pretzel. \
\textbf{Teacher}: Do you think it’s the same as a circle? \
\textbf{Student}: No way, it’s messier. I tried to undo it, but it stays twisted. \
\textbf{Teacher}: What if I showed you a different drawing, like this one? \
\textit{(Shows equivalent diagram)} \
\textbf{Student}: Whoa… it looks different, but the crossings are kind of the same. \
\textbf{Teacher}: Exactly! That’s the big question in knot theory: when are two knots the same?
\end{tcolorbox}

\begin{tcolorbox}[title=Classroom Snapshot 2: Tricoloring Frustration and Breakthrough, colback=gray!5, colframe=gray!75!black]
\textbf{Student}: I tried all the colors, but it doesn’t work. None of the rules match at all the crossings. \
\textbf{Teacher}: Interesting. What might that tell us about this knot? \
\textbf{Student}: (pauses) Maybe it’s not tricolorable like the trefoil? \
\textbf{Teacher}: That’s a solid hypothesis. Can you explain why? \
\textbf{Student}: Yeah! I tested every combo, and the rule about needing at least two colors at each crossing always fails.
\end{tcolorbox}

\begin{tcolorbox}[title=Classroom Snapshot 3: Reflection Moment, colback=gray!5, colframe=gray!75!black]
At the end of the second lesson, a student asked, “Is this math or is this art?”\footnote{This exchange is based on real conversations the author had with undergraduate students in the Business Analytics BSc program at Maastricht University.} \\
The teacher smiled and replied, “Maybe both. What do you think?” \\
“It’s like solving puzzles but with pictures”, the student replied. “It feels more creative than anything we’ve done so far.”
\end{tcolorbox}

\section{Discussion: Benefits and Challenges}

Knot theory offers more than just an engaging mathematical topic, it represents a pedagogical gateway to a richer, more interdisciplinary vision of education. As with any innovative approach, its classroom integration brings both significant advantages and notable challenges. In this section, we reflect on the broader educational benefits of knot theory through the lens of two influential frameworks, namely, Problem-Based Learning (PBL) and Global Citizenship Education (GCE), and consider their alignment with the VIBE principles introduced earlier. These pedagogical lenses illuminate knot theory's potential to enhance students' conceptual understanding, foster collaborative inquiry, and connect mathematics to global and cultural contexts. At the same time, we explore practical constraints such as curricular fit, teacher preparedness, and assessment.

\subsection*{Problem-Based Learning (PBL) in the Knot Theory Classroom}

As an institution that champions Problem-Based Learning (PBL), Maastricht University offers a pedagogical ethos deeply aligned with the goals of this paper. PBL is an instructional approach centered around active student inquiry, collaborative learning, and real-world problem solving. Rather than passively absorbing information, students in a PBL environment construct knowledge by working together on authentic tasks, guided by a tutor or facilitator.

Knot theory is particularly well-suited to PBL for several reasons:

\begin{enumerate}
\item \textbf{Authentic Problems}: Students can engage with open-ended questions such as Are these two knots equivalent?'', How can a knot be untangled?'', or ``What is the mathematical meaning of a knot invariant?'' These questions do not have a single correct answer and invite investigation.
\item \textbf{Collaborative Inquiry}: Learners work in groups to analyze knot diagrams, identify patterns, and explore equivalences. This promotes discussion, peer explanation, and shared meaning-making.
\item \textbf{Self-Directed Learning}: Students are encouraged to identify gaps in their understanding (e.g., needing to learn what a Reidemeister move is) and seek out resources collaboratively. The teacher takes the role of facilitator rather than lecturer.
\item \textbf{Integration of Skills}: PBL promotes skills beyond content knowledge, including communication, time management, reflective thinking, and critical analysis. In the context of knot theory, students must visualize, abstract, justify, and explain their reasoning to others.
\item \textbf{Alignment with VIBE Framework}: The core values of the VIBE framework (Visual Learning, Inquiry-Based Exploration, Braided Activities, and Embedded Contexts) resonate directly with the principles of PBL.
Visual tools serve as shared references for collaborative discussion; inquiry encourages open-ended questioning and personal engagement with knot concepts; braided activities foster peer learning and group-based construction of understanding; and embedded contexts help students anchor abstract ideas in real-world and interdisciplinary settings. Together, VIBE and PBL offer a synergistic blueprint for facilitating deep, active, and inclusive mathematical learning.
\end{enumerate}

By incorporating knot theory into a PBL classroom, educators can provide a powerful experience that blends advanced mathematical ideas with meaningful, student-driven exploration. This combination has the potential to transform mathematical learning into a deeply personal and collaborative endeavor.

\smallbreak 

\noindent\textbf{Historical and Interdisciplinary Contexts as PBL Catalysts.}
In a PBL environment, context is not a backdrop, it is a driver of curiosity. The rich history of knot theory provides fertile ground for meaningful problem scenarios. For instance, students might explore Lord Kelvin’s 19th-century hypothesis that atoms are knotted vortex tubes, which led Peter Guthrie Tait to construct early knot tables. While the physical theory was later debunked, the mathematical legacy remains. This episode offers an excellent prompt for inquiry: How can mistaken physical assumptions lead to fruitful mathematical frameworks?

Reidemeister’s formalization of knot equivalence through three elementary moves lends itself naturally to hands-on problem solving and collaborative exploration, core tenets of PBL. Students can be asked to determine whether two knots are equivalent using only these rules, creating opportunities for discovery, justification, and strategy comparison.

The modern mathematical significance of knots, such as the Jones polynomial linking topology, algebra, and quantum physics, offers an exciting challenge that need not remain out of reach for secondary students. While these invariants are often considered too advanced, we argue that a combinatorial and diagrammatic approach, adapted from Kauffman's skein-theoretic perspective \cite{LK1}, can make them accessible. By exploring how small local changes in a knot diagram (e.g., switching a crossing) affect an associated quantity, students can begin to grasp the recursive structure of polynomial invariants. With appropriate scaffolding, even high school learners can compute simple cases and interpret their results. This process empowers students to see themselves not just as consumers of mathematics, but as creators of it, capable of uncovering structure and inventing their own classification schemes. Such tasks blend abstraction with tangible exploration, making deep mathematics both meaningful and engaging.

Interdisciplinary connections also support rich, authentic problem contexts. DNA replication, for instance, leads to natural questions about how enzymes resolve entanglements, a perfect scenario to launch a biological modeling task using knots and links. Cultural motifs, from Celtic designs to Chinese and Andean traditions, allow students to explore how knot structures have encoded meaning across time and societies. Such explorations connect mathematics to identity, art, and global citizenship.

Thus, instead of treating history or applications as external anecdotes, this approach embeds them within the learning process. The goal is not only to solve problems, but to understand why these problems matter mathematically, culturally, and scientifically. In this way, the historical and interdisciplinary dimensions of knot theory become engines for inquiry and powerful levers for inclusive, story-driven mathematics education.

\subsection*{Global Citizenship Education and Broader Learning Outcomes}

Knot theory can also serve as a gateway to broader themes of Global Citizenship Education (GCE), which promotes values such as empathy, cultural awareness, and the ability to engage with complexity \cite{UNESCO2015, Nussbaum2006}. The study of knots, appearing in cultural artifacts, scientific models, and everyday practices around the world, naturally invites discussion about the universality and diversity of mathematical thinking.

By working on open problems, engaging in respectful dialogue, and drawing connections across disciplines and cultures, students practicing knot theory in the classroom are not only developing mathematical reasoning but also cultivating the skills of active global citizenship. These include critical thinking, cross-cultural communication, and a disposition toward collaboration, competencies that resonate deeply with 21st-century educational goals and align with the ethos of schools like Maastricht University.

\subsection*{Summary of Benefits}

The integration of knot theory into secondary mathematics education offers a wealth of pedagogical, cognitive, and social benefits. Drawing from the theoretical foundations outlined above and the practical examples showcased through the VIBE framework, we summarize the educational gains across multiple domains:

\begin{enumerate}
    
    \item \textbf{Advances Spatial and Visual Reasoning}:
Knot theory accommodates a wide range of learners through visual, tactile, and narrative entry points \cite{Duval1999, devlin}. Through manipulating ropes, drawing diagrams, and interpreting projections, students strengthen their ability to mentally rotate, deform, and compare complex structures. These skills are critical in fields ranging from engineering and architecture to biology and data visualization, yet often underdeveloped in conventional curricula.

    \item \textbf{Promotes Conceptual Depth and Abstract Reasoning}: Unlike rote exercises, knot theory tasks invite students to grapple with core mathematical ideas such as transformation, equivalence, and invariance. By distinguishing between physical deformations that preserve or alter a knot, learners engage with the essence of topology. These abstract concepts nurture flexible thinking and the ability to reason beyond symbols and formulas. 

    \item \textbf{Encourages Inquiry, Exploration, and Creativity}:
    The open-ended nature of knot theory problems, such as asking whether two knots are equivalent, creating new knot diagrams, or discovering invariants, supports a classroom culture of curiosity. Rather than seeking a single correct answer, students are invited to hypothesize, test, and refine their ideas. This mirrors authentic mathematical research and fosters creative problem-solving.

    \item \textbf{Enables Tactile and Multisensory Engagement}:
    The hands-on nature of knot activities (tying, manipulating, and drawing) provides a tangible entry point for students who may struggle with abstract or symbolic representations. Kinesthetic learners, in particular, benefit from experiencing mathematical ideas physically, which enhances memory, understanding, and confidence.

    \item \textbf{Supports Inclusive and Differentiated Instruction}:
    Knot theory is accessible to a wide range of learners. Visual thinkers, kinesthetic learners, and students who may be disengaged by algebraic abstraction can find meaningful entry points into mathematical thinking. Tasks can be easily scaffolded or extended, enabling differentiated instruction that challenges advanced students while supporting those who need more guidance.

    \item \textbf{Facilitates Interdisciplinary Connections}:
    Knot theory provides authentic bridges to other disciplines: DNA folding and protein entanglement in biology, woven design in visual art, encryption and routing in computer science, and naval applications in history and geography. These connections enrich the learning experience and illustrate mathematics as a vibrant, integrated part of human knowledge.

    \item \textbf{Aligns with Student-Centered Pedagogies}:
    The VIBE framework reflects contemporary educational priorities that value active learning, collaboration, and real-world relevance. Knot theory lessons model best practices in constructivist pedagogy and demonstrate how mathematics can be taught in dynamic, student-driven ways.

    \item \textbf{Encourages Cultural Awareness and Global Competence}:
    Knots are not only mathematical but also cultural artifacts. By exploring Celtic knots, Chinese decorative knots, or maritime traditions, students encounter mathematics through the lens of diverse cultural expressions. These discussions support Global Citizenship Education (GCE) goals, fostering respect for cultural diversity and appreciation of the universal language of form and pattern.

\end{enumerate}

To summarize the above, knot theory does more than introduce a new mathematical topic, it transforms how students experience and relate to mathematics. Through visualization, manipulation, collaboration, and contextualization, it cultivates a deeper, more personal connection to mathematical thinking.

\subsection*{Practical and Structural Challenges}

While knot theory offers compelling educational advantages, its implementation in high school classrooms faces several practical and systemic hurdles. These challenges span curriculum design, teacher readiness, assessment structures, and issues of equity. Addressing them is crucial for the sustainable integration of knot theory into broader educational practice.

\begin{enumerate}

    \item \textbf{Curricular Inflexibility}:
    Highly centralized education systems pose difficulties for the inclusion of non-standard topics \cite{Schmidt2006}. National and regional mathematics curricula are often tightly structured, with a strong emphasis on algebra, calculus, and statistics. These frameworks leave little room for exploratory or non-standard topics such as knot theory. Teachers may struggle to justify devoting class time to content not directly aligned with standardized learning outcomes. Integration may therefore require creative curricular embedding, such as interdisciplinary projects, enrichment modules, or extracurricular clubs.

    \item \textbf{Limited Teacher Preparation}:
    Topology and knot theory remain largely absent from secondary curricula, in part due to lack of teacher exposure \cite{murasugi}. Most secondary mathematics teachers receive minimal exposure to topology during their formal training. Knot theory, though accessible in spirit, can appear intimidating due to its association with advanced mathematics. Teachers may lack both the conceptual background and pedagogical strategies needed to confidently introduce knot-based lessons. Addressing this challenge requires targeted professional development opportunities, collaborative lesson planning, and the creation of accessible teaching guides. 

    \item \textbf{Assessment Misalignment}:
    Traditional assessment models prioritize procedural fluency, standardized formats, and closed-form answers. Knot theory, by contrast, emphasizes spatial reasoning, conjecture, and process-based learning. Without support for alternative assessment methods, such as visual journals, project-based tasks, or oral presentations, teachers may find it difficult to evaluate student understanding in ways that align with institutional expectations.

    \item \textbf{Resource Constraints}:
    Although many knot theory activities are low-cost, effective implementation benefits from visual aids, manipulatives, and digital tools that help students visualize and manipulate knots in two and three dimensions. Not all classrooms have access to these resources. Broader dissemination of open-access tools, printable materials, and simple physical models (e.g., string and cardboard templates) would help bridge this gap.

    \item \textbf{Perceptions of Relevance and Rigor}:
    Some educators and stakeholders may view knot theory as mathematically peripheral or recreational. Without a clear connection to assessment metrics or university preparation, there may be skepticism about its value. Advocating for knot theory requires articulating its benefits in terms of cognitive development, curricular enrichment, and alignment with broader educational goals such as inquiry-based learning and global competence.

    \item \textbf{Equity and Access Considerations}:
    Innovative topics like knot theory risk becoming exclusive to well-resourced schools or high-performing student groups unless deliberate efforts are made to ensure accessibility. Larger class sizes, limited planning time, or constrained curricular autonomy can inhibit implementation in underserved contexts. Ensuring equitable access requires both institutional support and scalable, adaptable teaching materials.

    \item \textbf{Managing Conceptual Complexity}:
    While knot theory begins with tangible and intuitive ideas, it can quickly lead to abstract or technically demanding content. Teachers must navigate a careful balance, offering depth and challenge without overwhelming students. Clear scaffolding, carefully sequenced activities, and a focus on intuition over formalism can help maintain accessibility.

\end{enumerate}

These challenges do not negate the pedagogical value of knot theory; rather, they highlight the need for strategic planning and institutional backing. With the right supports, professional development, adaptable resources, and curricular flexibility, knot theory can become a powerful tool for enriching mathematics education in diverse learning environments.

\subsection*{SWOT Analysis: Integrating Knot Theory through the VIBE Framework}

\begin{center}
\renewcommand{\arraystretch}{1.5}
\setlength{\tabcolsep}{8pt}
\begin{tabularx}{\textwidth}{|X|X|}
\hline
\rowcolor{gray!15}
\textbf{Strengths} & \textbf{Weaknesses} \\
\hline
Highly visual and tactile entry point into abstract mathematics & Requires time for teacher training and familiarization \\
\hline
Promotes inquiry, collaboration, and cross-disciplinary thinking & Limited space in standard curricula \\
\hline
Inclusive for diverse learners (visual, kinesthetic, creative) & Not included in standardized assessments \\
\hline
Supports 21st-century skills like spatial reasoning and creativity & May be perceived as recreational or non-essential by stakeholders \\
\hline
\rowcolor{gray!15}
\textbf{Opportunities} & \textbf{Threats} \\
\hline
Can be integrated into STEM/STEAM, math clubs, or interdisciplinary units & Uneven access to resources across schools \\
\hline
Opens doors to international, cultural, and artistic connections & Resistance to pedagogical innovation in rigid school systems \\
\hline
Sparks student-led research or creative projects & Risk of superficial treatment without proper scaffolding \\
\hline
Alignment with modern educational trends (PBL, UDL, GCE) & Competing priorities for curriculum time and policy constraints \\
\hline
\end{tabularx}
\end{center}

\subsection*{Moving Forward}

Despite the structural and practical obstacles outlined above, the integration of knot theory into secondary education represents a compelling opportunity to reimagine the teaching and learning of mathematics. It offers a shift in paradigm away from the narrow focus on speed, memorization, and procedural mastery, and toward a vision that values depth, creativity, and student-driven exploration.

Knot theory encourages learners to see mathematics as dynamic and visual, grounded in tangible experience and open-ended inquiry. Its emphasis on reasoning, transformation, and structure promotes a kind of mathematical thinking that is often underrepresented in traditional curricula. For many students, especially those who feel alienated by abstract symbols or formulaic exercises, knot theory can serve as a gateway into authentic mathematical engagement.

To realize this potential, several conditions must be cultivated. First, educators need access to adaptable resources, clear lesson plans, and opportunities for professional learning. Second, curricular structures must allow for exploratory content that connects across disciplines. Third, a broader cultural shift is needed, one that embraces mathematics not merely as a tool for calculation, but as a field of human creativity, beauty, and meaning.

With thoughtful design, strategic support, and a community of practice, knot theory can become more than a classroom novelty. It can serve as a catalyst for pedagogical renewal, helping students not only learn mathematics, but also see themselves as active participants in the mathematical world. The challenge is not whether knot theory belongs in schools, but how we can most effectively harness its potential to inspire the next generation of learners.

\subsection*{Future Research and Directions}

The integration of knot theory into secondary education opens the door to a broader set of questions, not only about pedagogy in mathematics, but also about curriculum design, educational equity, and interdisciplinary teaching more generally. To fully understand and optimize the impact of this approach, further research is needed.

One key avenue of investigation is the broader applicability of the VIBE framework. While developed here in the context of knot theory, the core principles, visual learning, inquiry-based exploration, braided (integrated and collaborative) activities, and embedded contextualization, may be transferable to other disciplines. For example:

\begin{itemize}
    \item Can the VIBE model enhance the teaching of abstract topics in physics, such as wave functions or quantum mechanics?
    \item Could VIBE be adapted to humanities subjects, fostering inquiry and collaboration in history, literature, or philosophy through visual narratives and thematic connections?
    \item What adjustments would be necessary to apply VIBE in courses with heavy symbolic or procedural content, such as algebra or chemistry?
\end{itemize}

In addition, future research could explore:

\begin{itemize}
    \item Longitudinal studies assessing the cognitive, affective, and social outcomes of students exposed to knot theory through the VIBE approach.
    \item Comparative studies between VIBE-based instruction and traditional pedagogies, focusing on motivation, retention, and conceptual understanding.
    \item Teacher training and professional development models that support the implementation of visually and topologically rich content in diverse classroom settings.
    \item The design and dissemination of digital tools, manipulatives, and cross-disciplinary resources that facilitate scalable adoption of knot-based curricula.
    \item The role of cultural and linguistic differences in students’ interpretations of visual mathematical concepts, and how VIBE might adapt to diverse educational contexts globally.
\end{itemize}

Ultimately, the future of this work lies in collaborative experimentation and interdisciplinary dialogue. By situating knot theory within a broader educational vision, we invite educators, researchers, and policymakers to rethink how we structure school mathematics, and how we might better align it with creativity, curiosity, and the human drive to understand complex systems through multiple lenses.

\section{Conclusion}

The integration of knot theory into secondary education represents a promising and transformative step toward more engaging, inclusive, and conceptually rich mathematics instruction. Through its visual, tactile, and inquiry-driven nature, knot theory bridges the abstract with the tangible, offering students meaningful entry points into higher-level mathematical thinking. As a pedagogical tool, it aligns naturally with constructivist learning principles, supports cross-disciplinary exploration, and nurtures skills that are essential for the 21st-century learner: reasoning, creativity, collaboration, and persistence.

The VIBE framework provides a cohesive and flexible model for teachers to implement knot theory in varied classroom settings. From the exploration of knot diagrams and Reidemeister moves to hands-on activities and real-world connections, the sample lessons outlined in this paper are intended to spark curiosity, foster conceptual understanding, and cultivate joy in mathematical discovery.

At the same time, this paper has emphasized the practical considerations for implementation: the importance of teacher support, adaptable assessment methods, inclusive practices, and resource availability. While challenges exist, they are not insurmountable. In fact, they offer opportunities for educational innovation and collaboration among educators, mathematicians, and curriculum developers.

Ultimately, by bringing knot theory into the high school curriculum, we invite students to think deeply, to see mathematics as a dynamic and creative endeavour, and to engage with one of the most beautiful ideas in topology. In doing so, we help redefine what mathematics education can be, not just a collection of rules and procedures, but a landscape of ideas, patterns, and human curiosity.

Future research may investigate how the VIBE framework can be generalized to other mathematical topics or disciplines, what supports are needed to scale implementation, and how student outcomes evolve over time with sustained exposure to such exploratory models. Knot theory, then, is not just a content choice, it is a case study in how reimagining mathematics education can help us build classrooms where all students are empowered to do mathematics with purpose and imagination.

\section*{Acknowledgments}

This paper was inspired during the preparation of my teaching portfolio for the BKO (University Teaching Qualification) at Maastricht University, The Netherlands. I am especially grateful to the team at the SBE Learning Academy for their invaluable support throughout the process. The reflective nature of the BKO experience played a key role in shaping the pedagogical perspectives developed in this work.

I would also like to thank Dr. Serdar T\"urkeli for encouraging me to approach academic responsibilities with a reflective and creative mindset. His advice to treat institutional tasks as opportunities for deeper inquiry helped reframe the BKO process as a meaningful exploration of pedagogy and practice, an approach that ultimately led to the development of this paper.

\end{document}